\RequirePackage{fix-cm}
\RequirePackage{fixltx2e}
\documentclass[oneside,english]{amsart}

\usepackage[latin9]{inputenc}
\usepackage{babel}
\usepackage{url}
\usepackage{amsthm}
\usepackage{amssymb}
\usepackage{esint}
\usepackage[unicode=true,pdfusetitle,
 bookmarks=true,bookmarksnumbered=false,bookmarksopen=false,
 breaklinks=false,pdfborder={0 0 1},backref=false,colorlinks=false]
 {hyperref}
 \makeatletter

\renewcommand{\theequation}{\theequation. \arabic{equation}}
\numberwithin{equation}{section}
\newtheorem{thm}{Theorem}[section]
\newtheorem{cor}[thm]{Corollary}

\newtheorem{prop}{Proposition}[section]

\usepackage[numbers,sort&compress]{natbib}

\def\squarebox#1{\hbox to #1{\hfill\vbox to #1{\vfill}}}
\def\qed{\hspace*{\fill}
         \vbox{\hrule\hbox{\vrule\squarebox{.667em}\vrule}\hrule}\smallskip}
\begin{document}\large
\title[Appell series $F_{3}$ over finite fields]
{\Large A finite field analogue for Appell series $F_{3}$}
\author{\small  Bing He}
\address{\small
College of Science, Northwest A\&F University,
   Yangling 712100, Shaanxi, People's Republic of China}
\email{yuhe001@foxmail.com; yuhelingyun@foxmail.com}




\begin{abstract}
\small In this paper we introduce a finite field analogue for the Appell series $F_{3}$ and
give some reduction formulae and  certain generating functions for this function over finite fields.
\end{abstract}
\maketitle
\section{Introduction}

Let $q$ be a power of a prime. $\mathbb{F}_{q}$ and $\widehat{\mathbb{F}}^{*}_{q}$ are denoted as the finite field of $q$ elements and  the group of multiplicative characters of $\mathbb{F}^{*}_{q}$ respectively. Then the domain of all characters $\chi$ of $\mathbb{F}^{*}_{q}$ can be extended to $\mathbb{F}_{q}$ by setting $\chi(0)=0$ for all characters. Let $\overline{\chi}$ and $\varepsilon$ denote  the inverse of $\chi$ and the trivial character respectively.  For more details about characters, please see \cite{BEW} and \cite [Chapter 8]{IR}.

Greene  in 1987 developed the theory of hypergeometric functions over finite fields and proved a number of transformation and summation identities for hypergeometric functions over finite fields which are analogues to those in the classical case \cite{Gr} (see \cite{B} for the definition of the hypergeometric functions).  Greene in \cite{Gr} introduced the notation
\begin{equation*}
{}_{2}F_1 \left(\begin{matrix}
A, B \\
C \end{matrix}
\bigg| x \right)^{G}=\varepsilon(x)\frac{BC(-1)}{q}\sum_{y}B(y)\overline{B}C(1-y)\overline{A}(1-xy)
\end{equation*}
for $A,B,C\in \widehat{\mathbb{F}}_{q}$ and $x\in \mathbb{F}_{q}$
and defined the finite field analogue of the binomial coefficient as
\begin{equation*}
  {A\choose B}^{G}=\frac{B(-1)}{q}J(A,\overline{B}),
\end{equation*}
where $J(\chi,\lambda)$ is the Jacobi sum given by $$J(\chi,\lambda)=\sum_{u}\chi(u)\lambda(1-u).$$ See \cite{FL, M, EG} for more information about the finite field analogue of the hypergeometric functions.

In this paper, for the sake of simplicity, we use the notation
\begin{equation*}
  {A\choose B}=q{A\choose B}^{G}=B(-1)J(A,\overline{B})
\end{equation*}
and define the finite field analogue of the classic Gauss hypergeometric function as
\begin{equation*}
{}_{2}F_1 \left(\begin{matrix}
A, B \\
C \end{matrix}
\bigg| x \right)=q\cdot{}_{2}F_1 \left(\begin{matrix}
A, B \\
C \end{matrix}
\bigg| x \right)^{G}=\varepsilon(x)BC(-1)\sum_{y}B(y)\overline{B}C(1-y)\overline{A}(1-xy).
\end{equation*}
Then
\begin{equation*}\label{e1-1}
  {}_{2}F_1 \left(\begin{matrix}
A, B \\
C \end{matrix}
\bigg| x \right)=\frac{1}{q-1}\sum_{\chi}{A\chi\choose \chi}{B\chi\choose C\chi}\chi(x)
\end{equation*}
for any $A,B,C\in \widehat{\mathbb{F}}_{q}$ and $x\in \mathbb{F}_{q}.$  Similarly, the finite field analogue of the generalized hypergeometric function for any $A_{0},A_{1},\cdots, A_{n},B_{1},\cdots,B_{n}\in \widehat{\mathbb{F}}_{q}$ and $x\in \mathbb{F}_{q}$ is defined by
\begin{equation*}
  {}_{n+1}F_n \left(\begin{matrix}
A_{0}, A_{1},\cdots, A_{n} \\
 B_{1},\cdots, B_{n} \end{matrix}
\bigg| x \right)=\frac{1}{q-1}\sum_{\chi}{A_{0}\chi\choose \chi}{A_{1}\chi\choose B_{1}\chi}\cdots{A_{n}\chi\choose B_{n}\chi}\chi(x).
\end{equation*}

In our notations one of Greene's theorems is as follows.
\begin{thm}\emph{(See \cite [Theorem 4.9]{Gr})}For any characters $A,B,C\in \widehat{\mathbb{F}_{q}},$  we have
\begin{equation}\label{e1-3}
{}_{2}F_1 \left(\begin{matrix}
A, B \\
 C \end{matrix}
\bigg| 1 \right)=A(-1){B\choose \overline{A}C}.
\end{equation}
\end{thm}

The results in the  following proposition  follows readily from some properties of Jacobi sums.
\begin{prop}\label{pp2-1} \emph{(See \cite [(2.6), (2.7), (2.8) and (2.13)]{Gr})} If $A,B\in \widehat{\mathbb{F}}_{q},$ then
\begin{align}
  {A\choose B}&={A\choose A\overline{B}},\label{f2}\\
  {A\choose B}&={B\overline{A}\choose B}B(-1)\label{f1},\\
  {A\choose B}&={\overline{B}\choose \overline{A}}AB(-1),\label{f3}\\
{\varepsilon\choose A}&=-A(-1)+(q-1)\delta(A). \label{f5}
\end{align}
where $\delta(\chi)$ is a function on characters given by
\begin{equation*}
\delta(\chi)=\left\{
               \begin{array}{ll}
                 1  & \hbox{if $\chi=\varepsilon$} \\
                 0  & \hbox{otherwise}
               \end{array}.
             \right.
\end{equation*}
\end{prop}

Among these interesting double hypergeometric functions in the field of hypergeometric functions, Appell's four functions may be the most important functions. Three of them are as follows.
\begin{align*}
F_{1}(a;b,b';c;x,y)&=\sum_{m,n\geq 0}\frac{(a)_{m+n}(b)_{m}(b')_{n}}{m!n!(c)_{m+n}}x^my^n,~|x|<1,~|y|<1,\\
F_{2}(a;b,b';c,c';x,y)&=\sum_{m,n\geq 0}\frac{(a)_{m+n}(b)_{m}(b')_{n}}{m!n!(c)_{m}(c')_{n}}x^my^n,~|x|+|y|<1,\\
F_{3}(a,a';b,b';c;x,y)&=\sum_{m,n\geq 0}\frac{(a)_{m}(a')_{n}(b)_{m}(b')_{n}}{m!n!(c)_{m+n}}x^my^n,~|x|<1,~|y|<1.
\end{align*}
See \cite{A, B, CA, S} for more detailed material about Appell's functions.

Inspired by Greene's work, Li \emph{et al} in \cite{LLM} gave a finite field analogue of the Appell series $F_{1}$ and
obtained some transformation and reduction formulas and the generating functions for the function over finite fields.  In that paper, the finite field analogue of the Appell series $F_{1}$ was given by
\begin{equation*}
F_{1}(A;B,B';C;x,y)=\varepsilon(xy)AC(-1)\sum_{u}A(u)\overline{A}C(1-u)\overline{B}(1-ux)\overline{B'}(1-uy).
\end{equation*}

Motivated by the work of Greene \cite{Gr} and Li \emph{et al} \cite{LLM}, the author \emph{et al} in \cite{HLZ} introduced a finite field analogue of  the Appell series  $F_{2}$  which is
\begin{align*}
  F_{2}(A;&B,B';C,C';x,y)\\
&=\varepsilon(xy)BB'CC'(-1) \sum_{u,v}B(u)B'(v)\overline{B}C(1-u)\overline{B'}C'(1-v)\overline{A}(1-ux-vy)
\end{align*}
and deduced certain transformation and reduction formulas and the generating functions for this function over finite fields.

In this paper we will give a finite field analogue for  the Appell series  $F_{3}.$ Since the Appell series  $F_{3}$ has the following  double integral representation  \cite [Chapter IX]{B}:
\begin{align*}
  F_{3}(a,a';b,b';c;x,y)&=\frac{\Gamma(c)}{\Gamma(b)\Gamma(b')\Gamma(c-b-b')}\\
&\cdot \int\int u^{b-1}v^{b'-1}(1-u-v)^{c-b-b'-1}(1-ux)^{-a}(1-vy)^{-a'}dudv,
\end{align*}
where the double integral is taken over the triangle region $\{(u,v)|u\geq 0, v\geq 0, u+v\leq 1\},$
we now give the finite field analogue of  $F_{3}$ in the form:
\begin{equation*}
F_{3}(A,A';B,B';C;x,y)=\varepsilon(xy)BB'(-1)\sum_{u,v}B(u)B'(v)C\overline{B}\overline{B'}(1-u-v)\overline{A}(1-ux)\overline{A'}(1-vy),
\end{equation*}
where $A,A',B,B',C\in \widehat{\mathbb{F}}_{q},~x,y\in \mathbb{F}_{q}$ and  each sum ranges over all the elements of $\mathbb{F}_{q}.$
In the above definition,  the factor $\frac{\Gamma(c)}{\Gamma(b)\Gamma(b')\Gamma(c-b-b')}$ is dropped to obtain simpler results. We choose the factor  $\varepsilon(xy)\cdot BB'(-1)$ to get a better expression in terms of binomial coefficients.

From the definition of $F_{3}(A,A';B,B';C;x,y)$  we know that
\begin{equation}\label{e1-2}
  F_{3}(A,A';B,B';C;x,y)=F_{3}(A',A;B',B;C;y,x)
\end{equation}
for any $A,B,B',C,C'\in \widehat{\mathbb{F}}_{q}$ and $x,y\in \mathbb{F}_{q}.$

The aim of this paper is to give several reduction formulas and  certain generating functions for the Appell series $F_{3}$ over finite fields. The fact that the Appell series $F_{3}$ does not  have  a single integral representation but has a  double one  leds us to giving  a finite field analogue for the Appell series $F_{3}$ which is more complicated than that for $F_{1}.$  Consequently,  the results on the reduction formulas and  the generating functions for the Appell series $F_{3}$ over finite fields are also quite complicated.

 We give two other expressions for $F_{3}(A,A';B,B';C;x,y)$ in the next section.  Several reduction formulae for $F_{3}(A,A';B,B';C;x,y)$ is given  in Section 3. The last section is devoted to deducing certain generating functions for $F_{3}(A,A';B,B';C;x,y).$
\section{Other expressions for $F_{3}(A,A';B,B';C;x,y)$}
In this section we give two other expressions for $F_{3}(A,A';B,B';C;x,y).$
\begin{thm}\label{t2-1}For any $A,A',B,B',C,\in \widehat{\mathbb{F}}_{q} $ and $x,y\in \mathbb{F}_{q},$ we have
\begin{align*}
  F_{3}(A,A';B,B';C;x,y)&=\frac{1}{(q-1)^2} \sum_{\chi, \lambda}{\overline{A}\choose \chi}{\overline{A'}\choose \lambda}{C\overline{B}\overline{B'}\choose C\overline{B'}\chi}{C\overline{B'}\chi\choose C\chi\lambda}\chi(x)\lambda(y)\\
&~+\varepsilon(y)B'(-1){\overline{A}\choose B'\overline{C}}B'\overline{C}(x)\overline{A'}(1-y),
\end{align*}
where each sum ranges over all multiplicative characters of $\mathbb{F}_{q}.$
\end{thm}
In order to prove Theorem \ref{t2-1} we need an auxiliary result.
\begin{prop}\emph{(Binomial theorem, see \cite [(2.5)]{Gr})} For any character $A\in \widehat{\mathbb{F}}_{q}$ and $x\in \mathbb{F}_{q},$ we have
\begin{equation*}
  A(1+x)=\delta(x)+\frac{1}{q-1}\sum_{\chi}{A\choose \chi}\chi(x),
\end{equation*}
where the sum ranges over all multiplicative characters of $\mathbb{F}_{q}$ and $\delta(x)$ is a function on $\mathbb{F}_{q}$ given by
\begin{equation*}
\delta(x)=\left\{
            \begin{array}{ll}
              1 & \hbox{if $x=0$} \\
              0 & \hbox{if $x\neq 0$}
            \end{array}.
          \right.
\end{equation*}
\end{prop}
We are now in the position to show Theorem \ref{t2-1}.\\
\emph{Proof of Theorem \ref{t2-1}.}
 When $v\neq 1,$ it is known from the binomial theorem that
\begin{equation*}
  C\overline{B}\overline{B'}\left(1-\frac{u}{1-v}\right)=\delta(u)+\frac{1}{q-1}\sum_{\mu}{C\overline{B}\overline{B'}\choose \mu}\mu(-u)\overline{\mu}(1-v).
\end{equation*}
Then
\begin{align}
\label{t2-1-1} C\overline{B}\overline{B'}(1-u-v)&=C\overline{B}\overline{B'}(1-v)C\overline{B}\overline{B'}\left(1-\frac{u}{1-v}\right)\\
 &=\delta(u)C\overline{B}\overline{B'}(1-v)+\frac{1}{q-1}\sum_{\mu}{C\overline{B}\overline{B'}\choose \mu}\mu(-u)C\overline{B}\overline{B'}\overline{\mu}(1-v).  \notag
\end{align}
It is easy to see from the binomial theorem that
\begin{align}
\overline{A}(1-ux)&=\delta(ux)+\frac{1}{q-1}\sum_{\chi}{\overline{A}\choose \chi}\chi(-ux), \label{t2-1-2}\\
\overline{A'}(1-vy)&=\delta(vy)+\frac{1}{q-1}\sum_{\lambda}{\overline{A'}\choose \lambda}\lambda(-vy). \label{t2-1-3}
\end{align}
Applying \eqref{t2-1-1}--\eqref{t2-1-3} and using that fact that $\delta(u)B(u)=\delta(ux)B(u)\varepsilon(xy)=\delta(vy)B'(v)\varepsilon(xy)=0,$  \cite [(1.15)]{Gr}  and \eqref{f2} yield
\begin{align*}
 &\varepsilon(xy)\sum_{u\in \mathbb{F}_{q},v\neq 1}B(u)B'(v)C\overline{B}\overline{B'}(1-u-v)\overline{A}(1-ux)\overline{A'}(1-vy)\\
 &=\frac{1}{(q-1)^3}\sum_{u\in \mathbb{F}_{q},v\neq 1}B(u)B'(v)\sum_{\mu}{C\overline{B}\overline{B'}\choose \mu}\mu(-u)C\overline{B}\overline{B'}\overline{\mu}(1-v)\\
 &~\cdot \sum_{\chi}{\overline{A}\choose \chi}\chi(-ux)\sum_{\lambda}{\overline{A'}\choose \lambda}\lambda(-vy)\\
 &=\frac{1}{(q-1)^3}\sum_{u\in \mathbb{F}_{q},v\in \mathbb{F}_{q}}B(u)B'(v)\sum_{\mu}{C\overline{B}\overline{B'}\choose \mu}\mu(-u)C\overline{B}\overline{B'}\overline{\mu}(1-v)\\
 &~\cdot \sum_{\chi}{\overline{A}\choose \chi}\chi(-ux)\sum_{\lambda}{\overline{A'}\choose \lambda}\lambda(-vy)\\
 &=\frac{1}{(q-1)^3}\sum_{\chi}{\overline{A}\choose \chi}\chi(-x)\sum_{\lambda}{\overline{A'}\choose \lambda}\lambda(-y)\sum_{\mu}{C\overline{B}\overline{B'}\choose \mu}\mu(-1)\\
 &~\cdot \sum_{u\in \mathbb{F}_{q}}B\chi\mu(u)\sum_{v\in \mathbb{F}_{q}}B'\lambda(v)C\overline{B}\overline{B'}\overline{\mu}(1-v)\\
 &=\frac{B(-1)}{(q-1)^2}\sum_{\chi}{\overline{A}\choose \chi}{C\overline{B}\overline{B'}\choose \overline{B}\overline{\chi}}\chi(x)\sum_{\lambda}{\overline{A'}\choose \lambda}\lambda(-y)\sum_{v}B'\lambda(v)C\overline{B'}\chi(1-v)\\
 &=\frac{BB'(-1)}{(q-1)^2} \sum_{\chi, \lambda}{\overline{A}\choose \chi}{\overline{A'}\choose \lambda}{C\overline{B}\overline{B'}\choose C\overline{B'}\chi}{C\overline{B'}\chi\choose C\chi\lambda}\chi(x)\lambda(y).
\end{align*}
On the other hand, by \eqref{t2-1-2}, the fact that $\varepsilon(xy)C\overline{B'}(u)\delta(ux)=0$  and \cite [(1.15)]{Gr}
\begin{align*}
&\varepsilon(xy)\sum_{u\in \mathbb{F}_{q},v= 1}B(u)B'(v)C\overline{B}\overline{B'}(1-u-v)\overline{A}(1-ux)\overline{A'}(1-vy)\\
&=\varepsilon(xy)C\overline{B}\overline{B'}(-1)\overline{A'}(1-y)\sum_{u\in \mathbb{F}_{q}}C\overline{B'}(u)\overline{A}(1-ux)\\
&=\frac{\varepsilon(y)C\overline{B}\overline{B'}(-1)\overline{A'}(1-y)}{q-1}\sum_{\chi}{\overline{A}\choose \chi}\chi(-x)\sum_{u\in \mathbb{F}_{q}}C\overline{B'}\chi(u),
\end{align*}
\begin{equation*}
=\varepsilon(y)B(-1){\overline{A}\choose B'\overline{C}}B'\overline{C}(x)\overline{A'}(1-y).
\end{equation*}
Therefore,
\begin{align*}
F_{3}(A,A';B,B';C;x,y)&=\varepsilon(xy)BB'(-1)\left(\sum_{u\in \mathbb{F}_{q},v\neq 1}+\sum_{u\in \mathbb{F}_{q},v=1}\right)\\
&=\frac{1}{(q-1)^2} \sum_{\chi, \lambda}{\overline{A}\choose \chi}{\overline{A'}\choose \lambda}{C\overline{B}\overline{B'}\choose C\overline{B'}\chi}{C\overline{B'}\chi\choose C\chi\lambda}\chi(x)\lambda(y)\\
&~+\varepsilon(y)B'(-1){\overline{A}\choose B'\overline{C}}B'\overline{C}(x)\overline{A'}(1-y).
\end{align*}
This completes the proof of Theorem \ref{t2-1}.\qed

\begin{cor}For any $A,A',B,B',C,\in \widehat{\mathbb{F}}_{q} $ and $x,y\in \mathbb{F}_{q},$ we have
\begin{align}
 F_{3}(A,A';B,B';C;x,1)&=B'C(-1){}_{3}F_2 \left(\begin{matrix}
A, B, \overline{A'}\overline{B'}C \\
 C\overline{B'}, C\overline{A'}  \end{matrix}
\bigg| x \right), \label{c2-1}\\
 F_{3}(A,A';B,B';C;1,y)&=BC(-1){}_{3}F_2 \left(\begin{matrix}
A', B', \overline{A}\overline{B}C \\
 C\overline{B}, C\overline{A}  \end{matrix}
\bigg| y \right).\label{c2-2}
\end{align}
\end{cor}
\noindent \emph{Proof.} We first show \eqref{c2-1}. It follows from Theorem \ref{t2-1}, \eqref{f1} and \eqref{e1-3} that
\begin{align*}
F_{3}(A,A';B,B';C;x,1)&=\frac{1}{(q-1)^2} \sum_{\chi, \lambda}{\overline{A}\choose \chi}{\overline{A'}\choose \lambda}{C\overline{B}\overline{B'}\choose C\overline{B'}\chi}{C\overline{B'}\chi\choose C\chi\lambda}\chi(x)\lambda(1)\\
&=\frac{C(-1)}{(q-1)^2} \sum_{\chi}{\overline{A}\choose \chi}{C\overline{B}\overline{B'}\choose C\overline{B'}\chi}\chi(-x)\sum_{\lambda}{A'\lambda\choose \lambda}{B'\lambda\choose C\chi\lambda}\lambda(1)\\
&=\frac{C(-1)}{q-1} \sum_{\chi}{\overline{A}\choose \chi}{C\overline{B}\overline{B'}\choose C\overline{B'}\chi}\chi(-x){}_{2}F_1 \left(\begin{matrix}
A', B' \\
C\chi \end{matrix}
\bigg| 1 \right)\\
&=\frac{A'C(-1)}{q-1} \sum_{\chi}{\overline{A}\choose \chi}{C\overline{B}\overline{B'}\choose C\overline{B'}\chi}{B'\choose C\overline{A'}\chi}\chi(-x)\\
&=\frac{B'C(-1)}{q-1} \sum_{\chi}{A\chi\choose \chi}{B\chi\choose C\overline{B'}\chi}{\overline{A'}\overline{B'}C\chi\choose C\overline{A'}\chi}\chi(x)\\
&=B'C(-1){}_{3}F_2 \left(\begin{matrix}
A, B, \overline{A'}\overline{B'}C \\
 C\overline{B'}, C\overline{A'}  \end{matrix}
\bigg| x \right),
\end{align*}
which proves \eqref{c2-1}. \eqref{c2-2} follows readily from \eqref{c2-1} and \eqref{e1-2}.\qed

The following theorem gives a third expression for $F_{3}(A,A';B,B';C;x,y).$
\begin{thm}\label{t2-2} For any $A,A',B,B',C,\in \widehat{\mathbb{F}}_{q} $ and $x,y\in \mathbb{F}_{q},$ we have
\begin{align*}
F_{3}(A,A';B,B';C;x,y)&=\frac{1}{(q-1)^2} \sum_{\chi, \lambda}{\overline{A}\choose \chi}{\overline{A'}\choose \lambda}{C\overline{B}\overline{B'}\choose C\chi\lambda}{\overline{B}\overline{B'}\overline{\chi}\overline{\lambda}\choose \overline{B'}\overline{\lambda}}\chi(x)\lambda(y)\\
&~+\frac{B'(-1)}{q-1} \sum_{\chi\lambda=\overline{B}\overline{B'}}{\overline{A}\choose \chi}{\overline{A'}\choose \lambda}\chi(x)\lambda(-y),
\end{align*}
where the sum in the second term of the right side ranges over the region $\chi, \lambda \in\widehat{\mathbb{F}}_{q}, \chi\lambda=\overline{B}\overline{B'}.$
\end{thm}
\noindent \emph{Proof.} It follows from the binomial theorem that
\begin{equation}\label{t2-2-1}
\sum_{\lambda}{\overline{A'}\choose \lambda}\lambda(-y)=(q-1)(\overline{A'}(1-y)-\delta(y))=(q-1)\varepsilon(y)\overline{A'}(1-y).
\end{equation}

According to \cite [(2.15)]{Gr}, we have
\begin{align*}
  {A\choose B}{C\choose A}={C\choose B}{C\overline{B} \choose A\overline{B}}-(q-1)B(-1)\delta(A)+(q-1)AB(-1)\delta(B\overline{C})
\end{align*}
for any $A,B,C\in \widehat{\mathbb{F}}_{q}.$  Then
\begin{align*}
 {C\overline{B}\overline{B'}\choose C\overline{B'}\chi}{C\overline{B'}\chi\choose C\chi\lambda}&={C\overline{B}\overline{B'}\choose C\chi\lambda}{\overline{B}\overline{B'}\overline{\chi}\overline{\lambda}\choose \overline{B'}\overline{\lambda}}-(q-1)C\chi\lambda(-1)\delta(C\overline{B'}\chi)\\
 &~+(q-1)B'\lambda(-1)\delta(BB'\chi\lambda).
\end{align*}
Using the above identity in Theorem  \ref{t2-1} and by \eqref{t2-2-1}, we get
\begin{align*}
  F_{3}(A,A';B,B';C;x,y)&=\frac{1}{(q-1)^2} \sum_{\chi, \lambda}{\overline{A}\choose \chi}{\overline{A'}\choose \lambda}{C\overline{B}\overline{B'}\choose C\chi\lambda}{\overline{B}\overline{B'}\overline{\chi}\overline{\lambda}\choose \overline{B'}\overline{\lambda}}\chi(x)\lambda(y)\\
  &~-\frac{C(-1)B'\overline{C}(-x)}{q-1} {\overline{A}\choose B'\overline{C}}\sum_{\lambda}{\overline{A'}\choose \lambda}\lambda(-y)\\
  &~+\frac{B'(-1)}{q-1} \sum_{\chi\lambda=\overline{B}\overline{B'}}{\overline{A}\choose \chi}{\overline{A'}\choose \lambda}\chi(x)\lambda(-y)\\
&~+\varepsilon(y)B'(-1){\overline{A}\choose B'\overline{C}}B'\overline{C}(x)\overline{A'}(1-y)\\
&=\frac{1}{(q-1)^2} \sum_{\chi, \lambda}{\overline{A}\choose \chi}{\overline{A'}\choose \lambda}{C\overline{B}\overline{B'}\choose C\chi\lambda}{\overline{B}\overline{B'}\overline{\chi}\overline{\lambda}\choose \overline{B'}\overline{\lambda}}\chi(x)\lambda(y)\\
&~+\frac{B'(-1)}{q-1} \sum_{\chi\lambda=\overline{B}\overline{B'}}{\overline{A}\choose \chi}{\overline{A'}\choose \lambda}\chi(x)\lambda(-y),
\end{align*}
which ends the proof of Theorem \ref{t2-2}. \qed

Actually, from Theorem \ref{t2-2} we can also deduce  \eqref{e1-2}, \eqref{c2-1} and \eqref{c2-2}.
\section{Reduction formulae}
In this section we give some reduction formulae for $F_{3}(A,A';B,B';C;
x,y).$
In order to derive these formulae we need some auxiliary results.
\begin{prop}\label{p2}\emph{(See \cite [Corollary 3.16 and Theorem 3.15]{Gr})} For any $A,B,C,D\in \widehat{\mathbb{F}}_{q}$ and $x \in \mathbb{F}_{q},$ we have
\begin{align}
{}_{2}F_1 \left(\begin{matrix}
 \varepsilon,B \\
C \end{matrix}
\bigg| x \right)&={B\choose C}\varepsilon(x)
-\overline{C}(x)\overline{B}C(1-x),\label{p31-4}\\
  {}_{2}F_1 \left(\begin{matrix}
A, \varepsilon \\
C \end{matrix}
\bigg| x \right)&={C\choose A}A(-1)\overline{C}(x)\overline{A}C(1-x)
-C(-1)\varepsilon(x)\label{p31-1}\\
&~~+(q-1)A(-1)\delta(1-x)\delta(\overline{A}C),\notag\\
{}_{2}F_1 \left(\begin{matrix}
A, B \\
A \end{matrix}
\bigg| x \right)&={B\choose A}\varepsilon(x)\overline{B}(1-x)-\overline{A}(-x)\label{p31-3}\\
&~~+(q-1)A(-1)\delta(1-x)\delta(B),\notag\\
{}_{3}F_2 \left(\begin{matrix}
A,B,C\\
D,B \end{matrix}
\bigg| x \right)&={C\overline{D}\choose B\overline{D}}{}_{2}F_1 \left(\begin{matrix}
A, C\\
D \end{matrix}
\bigg| x \right)-BD(-1)\overline{B}(x){A\overline{B}\choose \overline{B}}\label{p31-5}\\
&~~+(q-1)BD(-1)\delta(C\overline{D})\varepsilon(x)\overline{A}(1-x).\notag
\end{align}
\end{prop}

From the definition of $F_{3}(a,a';b,b';c;x,y)$ we know that
\begin{align*}
  F_{3}(a,0;b,b';c;x,y)&=F_{3}(a,a';b,0;c;x,y)={}_{2}F_1 \left(\begin{matrix}
a, b \\
c \end{matrix}
\bigg| x \right),\\
F_{3}(0,a';b,b';c;x,y)&=F_{3}(a,a';0,b';c;x,y)={}_{2}F_1 \left(\begin{matrix}
a, b' \\
c' \end{matrix}
\bigg| y \right).
\end{align*}
We now give finite field analogues of the above identities.
\begin{thm}\label{t31}Let $A,B,B', C,C'\in \widehat{\mathbb{F}}_{q}$ and $x,~y \in \mathbb{F}_{q}.$  If $y\neq 0,$ then
\begin{align}
  \label{t31-1} F_{3}(A,\varepsilon;B,B';C;x,y)&=CB'(-1) {B\overline{C}\choose \overline{B'}}{}_{2}F_1 \left(\begin{matrix}
A, B\\
C \end{matrix}
\bigg| x \right)-\overline{C}B'(x)\delta(y-1)B'(-1){\overline{A}\choose \overline{C}B'}\\
&~-B'(-1)\overline{C}(y)\overline{B'}C(1-y) {}_{2}F_1 \left(\begin{matrix}
 A,B \\
C\overline{B'} \end{matrix}
\bigg|-\frac{x(1-y)}{y} \right)\notag\\
&~+(q-1)\varepsilon(x)C(-1)\delta(B\overline{C})\overline{A}(1-x);\notag
\end{align}
if $x\neq 0,$ then
\begin{align}
\label{t31-2} F_{3}(\varepsilon, A';B,B';C;x,y)&=CB(-1) {B'\overline{C}\choose \overline{B}}{}_{2}F_1 \left(\begin{matrix}
A', B'\\
C \end{matrix}
\bigg| y \right)-\overline{C}B(y)\delta(x-1)B(-1){\overline{A'}\choose \overline{C}B}\\
&~-B(-1)\overline{C}(x)\overline{B}C(1-x) {}_{2}F_1 \left(\begin{matrix}
 A',B' \\
C\overline{B} \end{matrix}
\bigg|-\frac{y(1-x)}{x} \right)\notag\\
&~+(q-1)\varepsilon(y)C(-1)\delta(B'\overline{C})\overline{A'}(1-y). \notag
\end{align}
\end{thm}
\noindent{\it Proof.} We first prove \eqref{t31-1}. It is easily known from \eqref{p31-4} and \eqref{p31-5} that
\begin{align*}
  {}_{2}F_1 \left(\begin{matrix}
 \varepsilon,B' \\
C\chi \end{matrix}
\bigg| y \right)&={B'\choose C\chi}-\overline{C}\overline{\chi}(y)\overline{B'}C\chi(1-y),\\
{}_{3}F_2 \left(\begin{matrix}
A, C\overline{B'},B \\
C, C\overline{B'} \end{matrix}
\bigg| x \right)&={B\overline{C}\choose \overline{B'}}{}_{2}F_1 \left(\begin{matrix}
A, B\\
C \end{matrix}
\bigg| x \right)-B'(-1)\overline{C}B'(x){A\overline{C}B'\choose \overline{C}B'}\label{p22-5}\\
&~~+(q-1)B'(-1)\delta(B\overline{C})\varepsilon(x)\overline{A}(1-x).
\end{align*}
From Theorem \ref{t2-1}, the above two identities and \eqref{f1} we deduce that
\begin{align*}
F_{3}(A,\varepsilon;B,B';C;x,y)&=\frac{C(-1)}{(q-1)^2} \sum_{\chi}{\overline{A}\choose \chi}{C\overline{B}\overline{B'}\choose C\overline{B'}\chi}\chi(-x)
\sum_{\lambda}{\lambda\choose \lambda}{B'\lambda\choose C\chi\lambda}\lambda(y)\\
&~+B'(-1){\overline{A}\choose B'\overline{C}}B'\overline{C}(x)\varepsilon(1-y)\\
&=\frac{C(-1)}{q-1} \sum_{\chi}{\overline{A}\choose \chi}{C\overline{B}\overline{B'}\choose C\overline{B'}\chi}\chi(-x)
 {}_{2}F_1 \left(\begin{matrix}
 \varepsilon,B' \\
C\chi \end{matrix}
\bigg| y \right)\\
&~+B'(-1){\overline{A}\choose B'\overline{C}}B'\overline{C}(x)\varepsilon(1-y)\\
&=\frac{CB'(-1)}{q-1} \sum_{\chi}{A\chi\choose \chi}{B\chi\choose C\overline{B'}\chi}{C\overline{B'}\chi\choose C\chi}\chi(x)\\
&~-\frac{B'(-1)\overline{C}(y)\overline{B'}C(1-y)}{q-1} \sum_{\chi}{A\chi \choose \chi}{B\chi\choose C\overline{B'}\chi}\chi\left(-\frac{x(1-y)}{y}\right)\\
&~+B'(-1){\overline{A}\choose B'\overline{C}}B'\overline{C}(x)\varepsilon(1-y)\\
&=CB'(-1) {}_{3}F_2 \left(\begin{matrix}
A, C\overline{B'},B \\
C, C\overline{B'} \end{matrix}
\bigg| x \right)\\
&~-B'(-1)\overline{C}(y)\overline{B'}C(1-y) {}_{2}F_1 \left(\begin{matrix}
 A,B \\
C\overline{B'} \end{matrix}
\bigg|-\frac{x(1-y)}{y} \right)\\
&~+B'(-1){\overline{A}\choose B'\overline{C}}B'\overline{C}(x)\varepsilon(1-y)\\
&=CB'(-1) {B\overline{C}\choose \overline{B'}}{}_{2}F_1 \left(\begin{matrix}
A, B\\
C \end{matrix}
\bigg| x \right)-C(-1)\overline{C}B'(x){A\overline{C}B'\choose \overline{C}B'}\\
&~+(q-1)\varepsilon(x)C(-1)\delta(B\overline{C})\overline{A}(1-x)+B'(-1){\overline{A}\choose B'\overline{C}}B'\overline{C}(x)\varepsilon(1-y)\\
&~-B'(-1)\overline{C}(y)\overline{B'}C(1-y) {}_{2}F_1 \left(\begin{matrix}
 A,B \\
C\overline{B'} \end{matrix}
\bigg|-\frac{x(1-y)}{y} \right),
\end{align*}
which is equivalent to  \eqref{t31-1}. \eqref{t31-2} follows from \eqref{t31-1} and \eqref{e1-2}. This completes the proof of Theorem \ref{t31}.\qed
\begin{thm}\label{t3-2}Let $A,B,B', C,C'\in \widehat{\mathbb{F}}_{q}$ and $x,~y \in \mathbb{F}_{q}.$  If $y\neq 0,$ then
\begin{align*}
F_{3}(A,A';B,\varepsilon;C;x,y)&=-C(-1){}_{2}F_1 \left(\begin{matrix}
A, B\\
C \end{matrix}
\bigg| x\right)+{\overline{A}\choose A'\overline{C}}{\overline{B}C\choose A'}A'\overline{C}(x)\delta(1-y)\\
&~+A'(-1)\overline{C}(y)\overline{A'}C(1-y){A'B\overline{C}\choose A'}{}_{2}F_1 \left(\begin{matrix}
A, B\\
C\overline{A'} \end{matrix}
\bigg|-\frac{x(1-y)}{y}\right)\\
&~+(q-1)\varepsilon(x)\overline{C}(y)\overline{A'}C(1-y)\delta(A'B\overline{C})\overline{A}\left(1+\frac{x(1-y)}{y}\right);
\end{align*}
if $x\neq 0,$ then
\begin{align*}
F_{3}(A,A';\varepsilon,B';C;x,y)&=-C(-1){}_{2}F_1 \left(\begin{matrix}
A', B'\\
C \end{matrix}
\bigg| y\right)+{\overline{A'}\choose A\overline{C}}{\overline{B'}C\choose A}A\overline{C}(y)\delta(1-x)\\
&~+A(-1)\overline{C}(x)\overline{A}C(1-x){AB'\overline{C}\choose A}{}_{2}F_1 \left(\begin{matrix}
A', B'\\
C\overline{A} \end{matrix}
\bigg|-\frac{y(1-x)}{x}\right)\\
&~+(q-1)\varepsilon(y)\overline{C}(x)\overline{A}C(1-x)\delta(AB'\overline{C})\overline{A'}\left(1+\frac{y(1-x)}{x}\right).
\end{align*}
\end{thm}
\noindent{\it Proof.} It follows from \eqref{p31-1} and \eqref{p31-5} that
\begin{align*}
{}_{2}F_1 \left(\begin{matrix}
A', \varepsilon \\
C\chi \end{matrix}
\bigg| y \right)&={C\chi\choose A'}A'(-1)\overline{C}\overline{\chi}(y)\overline{A'} C\chi(1-y)-C\chi(-1)\\
&~+(q-1)A'(-1)\delta(1-y)\delta(\overline{A'}C\chi),\\
{}_{3}F_2 \left(\begin{matrix}
A,C, B\\
C\overline{A'}, C \end{matrix}
\bigg|-\frac{x(1-y)}{y}\right)&={A'B\overline{C}\choose A'}{}_{2}F_1 \left(\begin{matrix}
A, B\\
C\overline{A'} \end{matrix}
\bigg|-\frac{x(1-y)}{y}\right)\\
&~-A'(-1)C(y)\overline{C}(x(y-1)){A\overline{C}\choose \overline{C}}\\
&~+(q-1)A'(-1)\delta(A'B\overline{C})\varepsilon(x(1-y))\overline{A}\left(1+\frac{x(1-y)}{y}\right).
\end{align*}
We deduce from Theorem \ref{t2-1}, the above identities and \eqref{f1} that
\begin{align*}
  F_{3}(A,A';B,\varepsilon;C;x,y)&=\frac{1}{(q-1)^2} \sum_{\chi}{A\chi\choose \chi}{B\chi\choose C\chi}\chi(-x)\sum_{\lambda}{A'\lambda\choose \lambda}{\lambda\choose C\chi\lambda}\lambda(y)\\
&~+{\overline{A}\choose \overline{C}}\overline{C}(x)\overline{A'}(1-y)\\
&=\frac{1}{q-1} \sum_{\chi}{A\chi\choose \chi}{B\chi\choose C\chi}\chi(-x){}_{2}F_1 \left(\begin{matrix}
A', \varepsilon \\
C\chi \end{matrix}
\bigg| y \right)\\
&~+{\overline{A}\choose \overline{C}}\overline{C}(x)\overline{A'}(1-y)
\end{align*}
\begin{align*}
~~~~~~~~~~~~~~&=-\frac{C(-1)}{q-1} \sum_{\chi}{A\chi\choose \chi}{B\chi\choose C\chi}\chi(x)+{\overline{A}\choose \overline{C}}\overline{C}(x)\overline{A'}(1-y)\\
&~+\frac{A'(-1)\overline{C}(y)\overline{A'}C(1-y)}{q-1} \sum_{\chi}{A\chi\choose \chi}{B\chi\choose C\chi}{C\chi\choose C\overline{A'}\chi}\chi\left(-\frac{x(1-y)}{y}\right)\\
&~+ A'(-1){AA'\overline{C}\choose A'\overline{C}}{BA'\overline{C}\choose A'}A'\overline{C}(-x)\delta(1-y)\\
&=-C(-1){}_{2}F_1 \left(\begin{matrix}
A, B\\
C \end{matrix}
\bigg| x\right)+{\overline{A}\choose \overline{C}}\overline{C}(x)\overline{A'}(1-y)\\
&~+A'(-1)\overline{C}(y)\overline{A'}C(1-y){}_{3}F_2 \left(\begin{matrix}
A,C, B\\
C\overline{A'}, C \end{matrix}
\bigg|-\frac{x(1-y)}{y}\right)\\
&~+ A'(-1){AA'\overline{C}\choose A'\overline{C}}{BA'\overline{C}\choose A'}A'\overline{C}(-x)\delta(1-y)\\
&=-C(-1){}_{2}F_1 \left(\begin{matrix}
A, B\\
C \end{matrix}
\bigg| x\right)+{\overline{A}\choose A'\overline{C}}{\overline{B}C\choose A'}A'\overline{C}(x)\delta(1-y)\\
&~+A'(-1)\overline{C}(y)\overline{A'}C(1-y){A'B\overline{C}\choose A'}{}_{2}F_1 \left(\begin{matrix}
A, B\\
C\overline{A'} \end{matrix}
\bigg|-\frac{x(1-y)}{y}\right)\\
&~+(q-1)\varepsilon(x)\overline{C}(y)\overline{A'}C(1-y)\delta(A'B\overline{C})\overline{A}\left(1+\frac{x(1-y)}{y}\right),
\end{align*}
which proves the first identity. The second identity follows from the first identity and \eqref{e1-2}. This concludes the proof of Theorem \ref{t3-2}. \qed


\section{Generating functions}
In this section, we establish some generating functions for $F_{3}(A,A';B,B';C;x,y).$
\begin{thm}\label{t4-1}For any $A,A',B,B',C\in \widehat{\mathbb{F}}_{q}$ and $x,t \in \mathbb{F}^*_{q}\backslash\{1\},~y\in \mathbb{F}_{q},$ we have
\begin{align*}
&\frac{1}{q-1}\sum_{\theta}{A\theta\choose\theta}F_{3}(A\theta,A';B,B';C;x,y)\theta(t)\\
&=\overline{A}(1-t)F_{3}\left(A,A';B,B';C;\frac{x}{1-t},y\right)\\
&~-\overline{A}(-t) B(-1)\overline{C}(x)\overline{B}C(1-x){}_{2}F_1 \left(\begin{matrix}
A', B' \\
\overline{B}C \end{matrix}
\bigg| -\frac{y(1-x)}{x} \right).
\end{align*}
\end{thm}
\noindent{\it Proof.} It follows from \eqref{p31-3} and \eqref{f3} that
\begin{align*}
{}_{2}F_1 \left(\begin{matrix}
\overline{B'}\overline{\lambda}, B\overline{C}\overline{\lambda} \\
\overline{B'}\overline{\lambda} \end{matrix}
\bigg| x \right)&={B\overline{C}\overline{\lambda}\choose \overline{B'}\overline{\lambda}}\overline{B}C\lambda(1-x)-B'\lambda(-x)\\
&={B'\lambda\choose \overline{B}C\lambda}BB'C(-1)\overline{B}C\lambda(1-x)-B'\lambda(-x).
\end{align*}
Then, by \eqref{f1} and \cite [(2.11)]{Gr}
\begin{align}
&\sum_{\chi, \lambda}{\overline{A'}\choose \lambda}{C\overline{B}\overline{B'}\choose C\overline{B'}\chi}{C\overline{B'}\chi\choose C\chi\lambda}\chi(-x)\lambda(y)\\
&=CB'(-1)\sum_{\lambda}{A'\lambda\choose \lambda}\lambda(-y)\sum_{\chi}{B\chi\choose C\overline{B'}\chi}{C\overline{B'}\chi\choose C\chi\lambda}\chi(x) \notag\\
&=CB'(-1)\overline{C}(x)\sum_{\lambda}{A'\lambda\choose \lambda}\lambda\left(-\frac{y}{x}\right)\sum_{\chi}{B\overline{C}\overline{\lambda}\chi\choose \overline{B'}\overline{\lambda}\chi}{\overline{B'}\overline{\lambda}\chi\choose \chi}\chi(x) \notag\\
&=(q-1)CB'(-1)\overline{C}(x)\sum_{\lambda}{A'\lambda\choose \lambda}\lambda\left(-\frac{y}{x}\right){}_{2}F_1 \left(\begin{matrix}
\overline{B'}\overline{\lambda}, B\overline{C}\overline{\lambda} \\
\overline{B'}\overline{\lambda} \end{matrix}
\bigg| x \right) \notag\\
&=(q-1)B(-1)\overline{C}(x)\overline{B}C(1-x)\sum_{\lambda}{A'\lambda\choose \lambda}{B'\lambda\choose \overline{B}C\lambda}\lambda\left(-\frac{y(1-x)}{x}\right) \notag\\
&~-(q-1)\overline{C}(-x)B'(x)\sum_{\lambda}{A'\lambda\choose \lambda}\lambda(y) \notag\\
&=(q-1)^2 B(-1)\overline{C}(x)\overline{B}C(1-x){}_{2}F_1 \left(\begin{matrix}
A', B' \\
\overline{B}C \end{matrix}
\bigg| -\frac{y(1-x)}{x} \right) \notag\\
&~-(q-1)^2 \overline{C}(-x)B'(x)\varepsilon(y)\overline{A'}(1-y), \notag
\end{align}
where in the second step we have used the substitution $\chi\rightarrow \overline{C}\overline{\lambda}\chi.$
It is easily known from Theorem \ref{t2-1} that
\begin{align}
&\sum_{ \chi, \lambda}{\overline{A}\choose \chi}{\overline{A'}\choose \lambda}{C\overline{B}\overline{B'}\choose C\overline{B'}\chi}{C\overline{B'}\chi\choose C\chi\lambda}\chi\left(\frac{x}{1-t}\right)\lambda(y)\label{t41-2}\\
&=(q-1)^2\left(F_{3}\left(A,A';B,B';C;\frac{x}{1-t},y\right)-B'(-1)\varepsilon(y)B'\overline{C}(x)\overline{B'}C(1-t)\overline{A'}(1-y){\overline{A}\choose B'\overline{C}}\right). \notag
\end{align}
It follows from \eqref{p31-3} that
\begin{align}
\sum_{\theta}{A\theta\choose\theta}{A\chi\theta\choose A\theta}\theta(t)&=(q-1){}_{2}F_1 \left(\begin{matrix}
A, A\chi\\
A \end{matrix}
\bigg|t\right)\label{t41-3}\\
&=(q-1){A\chi\choose A}\overline{A}\overline{\chi}(1-t)-(q-1)\overline{A}(-t),\notag\\
\sum_{\theta}{A\theta\choose\theta}{AB'\overline{C}\theta\choose A\theta}\theta(t)&=(q-1){}_{2}F_1 \left(\begin{matrix}
A, AB'\overline{C}\\
A\end{matrix}
\bigg|t\right)\label{t41-4}\\
&=(q-1){AB'\overline{C}\choose A}\overline{A}\overline{B'}C(1-t)-(q-1)\overline{A}(-t).  \notag
\end{align}
By \eqref{f2}, \eqref{f3} and \eqref{t41-3}, we have
\begin{align}
&\sum_{ \chi, \lambda}{\overline{A'}\choose \lambda}{C\overline{B}\overline{B'}\choose C\overline{B'}\chi}{C\overline{B'}\chi\choose C\chi\lambda}\chi(-x)\lambda(y)\sum_{\theta}{A\theta\choose\theta}{A\chi\theta\choose A\theta}\theta(t) \label{t41-5}\\
&=(q-1)\overline{A}(1-t)\sum_{ \chi, \lambda}{\overline{A}\choose \chi}{\overline{A'}\choose \lambda}{C\overline{B}\overline{B'}\choose C\overline{B'}\chi}{C\overline{B'}\chi\choose C\chi\lambda}\chi\left(\frac{x}{1-t}\right)\lambda(y) \notag\\
&~-(q-1)\overline{A}(-t)\sum_{ \chi, \lambda}{\overline{A'}\choose \lambda}{C\overline{B}\overline{B'}\choose C\overline{B'}\chi}{C\overline{B'}\chi\choose C\chi\lambda}\chi(-x)\lambda(y). \notag
\end{align}
Substituting (4.1) and \eqref{t41-2} into \eqref{t41-5}, then applying Theorem \ref{t2-1}, \eqref{f2}, \eqref{f1}, \eqref{t41-4} and \eqref{t41-5} in the left side and cancelling some terms gives
\begin{align*}
 &\frac{1}{q-1}\sum_{\theta}{A\theta\choose\theta}F_{3}(A\theta,A';B,B';C;x,y)\theta(t)\\
 &=\frac{1}{(q-1)^3} \sum_{ \chi, \lambda}{\overline{A'}\choose \lambda}{C\overline{B}\overline{B'}\choose C\overline{B'}\chi}{C\overline{B'}\chi\choose C\chi\lambda}\chi(-x)\lambda(y)\sum_{\theta}{A\theta\choose\theta}{A\chi\theta\choose A\theta}\theta(t)\\
&~+\frac{1}{q-1}\varepsilon(y)C(-1)B'\overline{C}(x)\overline{A'}(1-y)\sum_{\theta}{A\theta\choose\theta}{AB'\overline{C}\theta\choose A\theta}\theta(t)\\
&=\overline{A}(1-t)F_{3}\left(A,A';B,B';C;\frac{x}{1-t},y\right)\\
&~-\overline{A}(-t) B(-1)\overline{C}(x)\overline{B}C(1-x){}_{2}F_1 \left(\begin{matrix}
A', B' \\
\overline{B}C \end{matrix}
\bigg| -\frac{y(1-x)}{x} \right).
\end{align*}
which is exactly the right side. This finishes the proof of Theorem \ref{t4-1}. \qed

From Theorem \ref{t4-1} and \eqref{e1-2} we can easily deduce another generating function for  $F_{3}(A,A';\\
B,B';C;x,y).$
\begin{thm}For any $A,A',B,B',C\in \widehat{\mathbb{F}}_{q}$ and $y,t \in \mathbb{F}^*_{q}\backslash\{1\},~x\in \mathbb{F}_{q},$ we have
\begin{align*}
&\frac{1}{q-1}\sum_{\theta}{A'\theta\choose\theta}F_{3}(A,A'\theta;B,B';C;x,y)\theta(t)\\
&=\overline{A'}(1-t)F_{3}\left(A,A';B,B';C; x, \frac{y}{1-t}\right)\\
&~-\overline{A'}(-t) B'(-1)\overline{C}(y)\overline{B'}C(1-y){}_{2}F_1 \left(\begin{matrix}
A, B \\
\overline{B'}C \end{matrix}
\bigg| -\frac{x(1-y)}{y} \right).
\end{align*}
\end{thm}
We also establish two other generating functions for $F_{3}(A,A'; B,B';C;x,y).$
\begin{thm}\label{t4-3} For any $A,A',B,B',C\in \widehat{\mathbb{F}}_{q}$ and $x,t \in \mathbb{F}^*_{q}\backslash\{1\},~y\in \mathbb{F}_{q},$ we have
\begin{align*}
&\frac{1}{q-1}\sum_{\theta}{BB'\overline{C}\theta\choose\theta}F_{3}(A,A';B\theta,B';C;x,y)\theta(t)\\
&=\overline{B}(1-t)F_{3}\left(A,A';B,B';C;\frac{x}{1-t},y\right)+\overline{B}\overline{B'}C(-t)F_{3}(A,A';C\overline{B'},B';C;x,y)\\
&~-\overline{B'}(y)\overline{B}(-t)\overline{C}B'\left(-\frac{x}{t}\right){\overline{A}\choose \overline{C}B'}{\overline{A'}\choose \overline{B'}}.
\end{align*}
\end{thm}
\noindent{\it Proof.} It is easily known from Theorem \ref{t2-1} that
\begin{align*}
F_{3}(A,A';C\overline{B'},B';C;x,y)&=\frac{1}{(q-1)^2}\sum_{\chi, \lambda}{\overline{A}\choose \chi}{\overline{A'}\choose \lambda}{\varepsilon\choose C\overline{B'}\chi}{C\overline{B'}\chi\choose C\chi\lambda}\chi(x)\lambda(y)\\
&~+\varepsilon(y)B'(-1){\overline{A}\choose B'\overline{C}}B'\overline{C}(x)\overline{A'}(1-y).
\end{align*}
Then
\begin{align}
&\sum_{\chi, \lambda}{\overline{A}\choose \chi}{\overline{A'}\choose \lambda}{\varepsilon\choose C\overline{B'}\chi}{C\overline{B'}\chi\choose C\chi\lambda}\chi(x)\lambda(y)\label{t43-1}\\
&=(q-1)^2\left(F_{3}(A,A';C\overline{B'},B';C;x,y)-\varepsilon(y)B'(-1){\overline{A}\choose B'\overline{C}}B'\overline{C}(x)\overline{A'}(1-y)\right). \notag
\end{align}
We see from \eqref{f1}, \eqref{f5}  and \cite [(2.10)]{Gr} that
\begin{align*}
\sum_{\lambda}{\overline{A'}\choose \lambda}{\varepsilon\choose B'\lambda}\lambda(y)&=-B'(-1)\sum_{\lambda}{A'\lambda \choose \lambda}\lambda(y)+(q-1){\overline{A'}\choose \overline{B'}}\overline{B'}(y)\\
&=-(q-1)B'(-1)\varepsilon(y)\overline{A'}(1-y)+(q-1){\overline{A'}\choose \overline{B'}}\overline{B'}(y).
\end{align*}
Then, by \eqref{f5},
\begin{align}
\label{t43-2} &\sum_{\chi, \lambda}{\overline{A}\choose \chi}{\overline{A'}\choose \lambda}{\varepsilon\choose C\overline{B'}\chi}{C\overline{B'}\chi\choose C\chi\lambda}\chi(x)\lambda(y)\\
&=-CB'(-1)\sum_{\chi, \lambda}{\overline{A}\choose \chi}{\overline{A'}\choose \lambda}{C\overline{B'}\chi\choose C\chi\lambda}\chi(-x)\lambda(y)+(q-1)\overline{C}B'(x){\overline{A}\choose \overline{C}B'}\sum_{\lambda}{\overline{A'}\choose \lambda}{\varepsilon\choose B'\lambda}\lambda(y)\notag\\
&=-CB'(-1)\sum_{\chi, \lambda}{\overline{A}\choose \chi}{\overline{A'}\choose \lambda}{C\overline{B'}\chi\choose C\chi\lambda}\chi(-x)\lambda(y)-(q-1)^2B'(-1)\overline{C}B'(x){\overline{A}\choose \overline{C}B'}\varepsilon(y)\overline{A'}(1-y)  \notag\\
&~+(q-1)^2\overline{C}B'(x){\overline{A}\choose \overline{C}B'}{\overline{A'}\choose \overline{B'}}\overline{B'}(y).\notag
\end{align}
From \eqref{t43-1} and \eqref{t43-2} we deduce that
\begin{align}
&\sum_{\chi, \lambda}{\overline{A}\choose \chi}{\overline{A'}\choose \lambda}{C\overline{B'}\chi\choose C\chi\lambda}\chi(-x)\lambda(y) \label{t43-3}\\
&=(q-1)^2\left(\overline{C}B'(-x){\overline{A}\choose \overline{C}B'}{\overline{A'}\choose \overline{B'}}\overline{B'}(y)-CB'(-1)F_{3}(A,A';C\overline{B'},B';C;x,y)\right).\notag
\end{align}
It follows from \eqref{p31-3} that
\begin{align*}
\sum_{\theta}{BB'\overline{C}\theta\choose\theta}{B\chi\theta\choose BB'\overline{C}\theta} \theta(t)&=(q-1){}_{2}F_1 \left(\begin{matrix}
BB'\overline{C}, B\chi \\
BB'\overline{C} \end{matrix}
\bigg|t\right)\\
&=(q-1){B\chi\choose BB'\overline{C}}\overline{B}\overline{\chi}(1-t)-(q-1)\overline{B}\overline{B'}C(-t).
\end{align*}
Then
\begin{align}
 &\sum_{\chi, \lambda}{\overline{A}\choose \chi}{\overline{A'}\choose \lambda}{C\overline{B'}\chi\choose C\chi\lambda}\chi(-x)\lambda(y)\sum_{\theta}{BB'\overline{C}\theta\choose\theta}{B\chi\theta\choose BB'\overline{C}\theta} \theta(t)\label{t43-4}\\
 &=(q-1)CB'(-1)\overline{B}(1-t)\sum_{\chi, \lambda}{\overline{A}\choose \chi}{\overline{A'}\choose \lambda}{C\overline{B}\overline{B'}\choose C\overline{B'}\chi}{C\overline{B'}\chi\choose C\chi\lambda}\chi\left(\frac{x}{1-t}\right)\lambda(y)\notag\\
 &~-(q-1)\overline{B}\overline{B'}C(-t)\sum_{\chi, \lambda}{\overline{A}\choose \chi}{\overline{A'}\choose \lambda}{C\overline{B'}\chi\choose C\chi\lambda}\chi(-x)\lambda(y). \notag
\end{align}
Substituting \eqref{t41-2} and \eqref{t43-3} into \eqref{t43-4}, applying \eqref{f2}, \eqref{f1}, \eqref{t43-4} and \cite [(2.10)]{Gr} in the left side  and cancelling some terms yields
\begin{align*}
  &\frac{1}{q-1}\sum_{\theta}{BB'\overline{C}\theta\choose\theta}F_{3}(A,A';B\theta,B';C;x,y)\theta(t)\\
  &=\frac{CB'(-1)}{(q-1)^3} \sum_{\chi, \lambda}{\overline{A}\choose \chi}{\overline{A'}\choose \lambda}{C\overline{B'}\chi\choose C\chi\lambda}\chi(-x)\lambda(y)\sum_{\theta}{BB'\overline{C}\theta\choose\theta}{B\chi\theta\choose BB'\overline{C}\theta} \theta(t)\\
&~+\frac{1}{q-1}\varepsilon(y)B'(-1){\overline{A}\choose B'\overline{C}}B'\overline{C}(x)\overline{A'}(1-y)\sum_{\theta}{BB'\overline{C}\theta\choose\theta} \theta(t)\\
&=\overline{B}(1-t)F_{3}\left(A,A';B,B';C;\frac{x}{1-t},y\right)\\
&~-\overline{B}(-t)\overline{B'}C(t)\left(\overline{C}B'(-x){\overline{A}\choose \overline{C}B'}{\overline{A'}\choose \overline{B'}}\overline{B'}(y)-CB'(-1)F_{3}(A,A';C\overline{B'},B';C;x,y)\right),
\end{align*}
which equals the right side. The proof of Theorem \ref{t4-3} is completed. \qed

From Theorem \ref{t4-3} and \eqref{e1-2} we can also deduce another generating function for  $F_{3}(A,A';\\
B,B';C;x,y).$
\begin{thm}For any $A,A',B,B',C\in \widehat{\mathbb{F}}_{q}$ and $y,t \in \mathbb{F}^*_{q}\backslash\{1\},~x\in \mathbb{F}_{q},$ we have
\begin{align*}
&\frac{1}{q-1}\sum_{\theta}{BB'\overline{C}\theta\choose\theta}F_{3}(A,A';B,B'\theta;C;x,y)\theta(t)\\
&=\overline{B'}(1-t)F_{3}\left(A,A';B,B';C;x, \frac{y}{1-t}\right)+\overline{B}\overline{B'}C(-t)F_{3}(A,A';B, C\overline{B};C;x,y)\\
&~-\overline{B}(y)\overline{B'}(-t)\overline{C}B\left(-\frac{y}{t}\right){\overline{A'}\choose \overline{C}B}{\overline{A}\choose \overline{B}}.
\end{align*}
\end{thm}
\section*{Acknowledgement}
This work was  supported by the Natural Science Basic Research Plan in Shaanxi Province of China (No. 2017JQ1001), the Initial Foundation for
Scientific Research  of Northwest A\&F University (No. 2452015321)  and   the Fundamental Research Funds for the Central Universities (No. 2452017170).

\end{document}